\documentclass[12pt,twoside,a4paper]{amsart}
\usepackage{amsmath,amssymb}
\usepackage{mathrsfs} 
\usepackage{mathtools}
\usepackage{cite}
\usepackage{graphicx}
\usepackage{txfonts}
\usepackage{mathrsfs}
\usepackage{empheq}
\newcommand{\ds}{\displaystyle}

\newcommand{\ka}{J.Kaczorowski}
\newcommand{\m}{H.L.Montgomery}

\newcommand{\se}{$\mathcal{S}$}

\newcommand{\tf}{f \otimes \chi}
\newcommand{\Ll}{\mathcal{L}}
\newcommand{\lL}{\mathscr{L}}
\newcommand{\la}[2]{\lambda_#1 \left( #2 \right)}

\allowdisplaybreaks[4]
\theoremstyle{definition}
\newtheorem{definition}{Definition}[subsection]
\newtheorem{remark}[definition]{Remark}

\theoremstyle{theorem}
\newtheorem{theorem}{Theorem}[subsection]
\newtheorem{lemma}[theorem]{Lemma}

\numberwithin{equation}{subsection}
\mathtoolsset{showonlyrefs=true}
\mathtoolsset{showonlyrefs=false}

\begin{document}

\title[A relation to a remainder terms \ldots]{A relation to a remainder terms in an asymptotic formula for the associated Euler totient function}
\author{Hideto Iwata}
\address{ National Institute of Technology, Gunma College, Japan.
}
\email{iwata@gunma-ct.ac.jp}
\date{\today}
\maketitle

\textbf{Abstract.} For any positive integer $n$, {\m} proved a relation for error terms in asymptotic  formulas for the Euler totient function $\varphi(n)$. {\ka} defined the associated Euler totient function which generalizes $\varphi(n)$ and obtained an asymptotic formula for it. In this paper, we prove a relation on error terms similar to {\m}'s result for a certain special case of the associated Euler totient function. 

\section{Introduction}

\subsection{Previous research and the Main Theorem}
Let $\varphi(n)$ be the Euler totient function and let
\begin{equation}
   E(x) = \sum_{n \leq x} \varphi(n) - \frac{3}{\pi^2}x^2   \label{eq:1.1.1} 
\end{equation}be the associated error term. Also, let
 
\begin{equation}
   E_0 (x) = \sum_{n \leq x} \frac{\varphi(n)}{n} - \frac{6}{\pi^2}x.   \label{eq:1.1.2} 
\end{equation}{\m} formulated the relation between these error terms in the following form.
\begin{theorem}[Theorem 1 in ~\cite{Mon}]
For $x \geq 2$,
\begin{equation}
   E_0 (x) = E(x)/x + O(\exp(-c\sqrt{\log x})).  \label{eq:1.1.3}
\end{equation}Here, $c$ is a positive absolute constant.
\end{theorem}The aim of this paper is to obtain a result similar to Theorem 1.1.1 for a function $F$ given in \eqref{eq:2.1.1} later whose Euler degree is equal to $2$. We consider the $L$-function attached to a Hecke-Maass form $f$ in Section 2.4 which is of the form $F$.The following relation similar to Theorem 1.1.1 holds for $E(x.,L)$ and $E_0 (x,L)$, where $E(x,L)$ and $E_0 (x,L)$ are error terms \eqref{eq:3.1.7} and \eqref{eq:3.1.8} given in Section 3.1.
\begin{theorem}
For $x \geq 2$,
\begin{equation}
    E_0 (x,L) = x^{-1}E(x,L) + O(\exp(-c_1 \sqrt{\log x})).   \label{eq:.1.1.4}
\end{equation}Here $c_1$ is a positive absolute constant.
\end{theorem}

\section{Backgrounds}
\subsection{The associated Euler totient function}
{\ka} defined the associated Euler totient function for a class of generalized $L$-functions including the Riemann zeta function, Dirichlet $L$-functions and obtained an asymptotic formula (see~\cite{Kac}) : By a polynomial Euler product we mean a function $F(s)$ of a complex variable $s=\sigma+it$ which for $\sigma>1$ is defined by the product of the form
\begin{equation}
   F(s) = \prod_{p} F_p (s) = \prod_{p}\prod_{j=1}^d \left( 1-\frac{\alpha_j (p)}{ p^s } \right)^{-1},   \label{eq:2.1.1}
\end{equation}where $p$ runs over primes and $|\alpha_j (p)| \leq 1$ for all $p$ and $1\leq j\leq d$. We assume that $d$ is chosen as small as possible, i.e. that there exists at least one prime number $p_0$ such that
\[ \ds \prod_{j=1}^d \alpha_j ( p_0 ) \neq 0. \] Then $d$ is called the \textit{Euler degree} of $F$. Note that the $L$-functions from number theory including the Riemann zeta function, Dirichlet $L$-functions, Dedekind zeta and Hecke $L$-functions of algebraic number fields, as well as the (normalized) $L$-functions of holomorphic modular forms and, conjecturally, many other $L$-functions are polynomial Euler products. For $F$ in \eqref{eq:2.1.1} we define \textit{the associated Euler totient function} as follows :
\begin{equation}
   \varphi(n,F) = n\prod_{p|n} F_p (1)^{-1} \quad (n \in \mathbb{N}).    \label{eq:2.1.2}
\end{equation} Let
\begin{gather}
   \gamma(p) = p\left( 1-\frac{1}{F_p (1)} \right),   \label{eq:2.1.3}
   \\
   C(F) = \frac{1}{2}\prod_p \left( 1 - \frac{\gamma(p)}{p^2} \right),   \label{eq:2.1.4}
\end{gather}and
\begin{equation}
   \alpha(n) = \mu(n)\prod_{p|n}\gamma(p).   \label{eq:2.1.5}
\end{equation}By \eqref{eq:2.1.1} and \eqref{eq:2.1.2}, we see that the Euler totient function corresponds to the case when $F$ is the Riemann zeta function $\zeta(s)$ whose Euler degree is equal to $1$. J.Kaczorowski obtained the following an asymptotic formula for \eqref{eq:2.1.2}.
\begin{theorem}[Theorem 1.1 in ~\cite{Kac}]
For a polynomial Euler product $F$ of degree $d$ and $x \geq 1$ we have
\begin{equation}
   \sum_{n \leq x} \varphi(n,F) = C(F)x^2 + O(x(\log 2x)^d).   \label{eq:2.1.6}
\end{equation}
\end{theorem}
\begin{remark}[p33 in ~\cite{Kac} ]Let us observe that for every positive $\epsilon$, $\alpha(n) \ll n^\epsilon$ in \eqref{eq:2.1.5}. Hence the series
\begin{equation}
   \sum_{n=1}^\infty \frac{\alpha(n)}{n^s}   \label{eq:2.1.7}
\end{equation}converges absolutely for $\sigma>1$.Since $\alpha(n)$ is multiplicative, we have
\begin{equation}
   \sum_{n=1}^\infty \frac{\alpha(n)}{n^2} = 2C(F).   \label{eq:2.1.8}
\end{equation}
\end{remark}

Now we provide the definition of the \textit{Selberg class} {\se} for our later purpose as follows (see ~\cite{P}) : $F \in$ {\se} if
\begin{enumerate}
\item[(i)](\textit{ordinary Dirichlet series}) $\displaystyle F(s) = \sum_{n=1}^\infty a_F (n)n^{-s}$, absolutely convergent for $\sigma > 1$; 
\item[(ii)](\textit{analytic continuation}) there exists an integer $m\geq0$ such that $(s-1)^m \cdot F(s)$ is an entire function of finite order; 
\item[(iii)](\textit{functional equation}) $F(s)$ satisfies a functional equation of type $\Phi(s) = \omega\overline{\Phi(1-\overline{s})}$, where
                                                           \begin{equation} 
                                                              \Phi(s) = Q^s \prod_{j=1}^r \Gamma(\lambda_j s + \mu_j)F(s) = \gamma(s)F(s),   \label{eq:2.1.9}
                                                           \end{equation} say, with $r\geq0, Q>0, \lambda_j >0$, Re\hspace{0.001cm} $\mu_j \geq 0$ and $|\omega| = 1$;
\item[(iv)](\textit{Ramanujan conjecture}) for every $\epsilon>0, a_F (n) \ll n^\epsilon$.  
\item[(v)](\textit{Euler product}) $\displaystyle F(s) = \prod_{p} \exp \left( \sum_{\ell=0}^\infty  \frac{b_F ( p^\ell )}{p^{\ell s}} \right)$, where $b_F (n) = 0$ unless $n=p^m$ with $m\geq1$, and $b_F (n) \ll n^\vartheta$ for some 
                                  $\vartheta<\frac{1}{2}$.
\end{enumerate}Note that we understand that an empty product is equal to $1$. 
\begin{theorem}[Lemma 2.2 in ~\cite{Kac}]
The series
\begin{equation}
   \sum_{n=1}^\infty \frac{\varphi(n,F)}{n^s}   \label{eq:2.1.10}
\end{equation}converges absolutely for $\sigma>2$ and in this half-plane we have
\begin{equation}
   \sum_{n=1}^\infty \frac{\varphi(n,F)}{n^s} = \zeta(s-1)\sum_{n=1}^\infty \frac{\alpha(n)}{n^s}.  \label{eq:2.1.11}
\end{equation}In particular,
\begin{equation}
   \varphi(n,F) = n\sum_{m|n}\frac{\alpha(m)}{m}.   \label{eq:2.1.12}
\end{equation}
\end{theorem}

\begin{theorem}[Lemma 2.3 in ~\cite{Kac}]
For $\sigma>1$ we have
\begin{equation}
   \sum_{n=1}^\infty \frac{\alpha(n)}{n^s} = \frac{H(s)}{F(s)},   \label{eq:2.1.13}
\end{equation}where 
\begin{equation}
   H(s) = \sum_{n=1}^\infty \frac{h(n)}{n^s}   \label{eq:2.1.14}
\end{equation}converges absolutely for $\sigma>1/2$. Moreover, as $n$ runs over square-free positive integers we have
\begin{equation}
   h(n) \ll \frac{1}{n}\exp\left( c^\prime \frac{\log n}{\log \log (n+2)} \right),   \label{eq:2.1.15}
\end{equation}where $c^\prime$ is a positive constant which may depend on $F$ and other parameters. In particular for such $n$, $h(n)$ is bounded.
\end{theorem}
\begin{remark}
We do not know a behavior of the function $H(s)$ for $\sigma \leq \frac{1}{2}$.  More details of $H(s)$ are written in~\cite{iwa}. If we take $\zeta(s)$ as $F$, $\alpha(n) = \mu(n)$ by \eqref{eq:2.1.5}. Hence, $H(s) = 1$ for all $s$ with $\sigma > 1$ by the Dirichlet series expansion for $\zeta(s)^{-1}$. 
\end{remark}

\subsection{Maass forms}
We prepare a brief account of the theory of Maass forms in Section 2.2-4 (see ~\cite{fg}) : Let $k$ be an integer, $q$ be a positive integer, and $\chi$ be a Dirichlet character modulo $q$ that satisfies the consistency condition $\chi(-1) = (-1)^k$. Such a character gives rise to a character of the Hecke congruence group 
\begin{equation*}
\Gamma_0 (q) = \left\{ \begin{pmatrix}
a & b \\
c & d \\
\end{pmatrix}
\in \text{SL}(2,\mathbb{Z}) \hspace{0.1cm} ; \hspace{0.1cm} c \equiv 0 \pmod q
 \right\}
\end{equation*}by $\chi(\gamma) = \chi(d)$ for $
\gamma = \begin{pmatrix}
a & b \\
c & d \\
\end{pmatrix}
\in \Gamma_0 (q)$. For $z \in \mathbb{H}$, the upper half plane, we set
\[ j_\gamma (z) := (cz+d)|cz+d|^{-1} = e^{i\arg(cz+d)}. \] A function $f : \mathbb{H} \longrightarrow \mathbb{C}$ that satisfies the condition
\[ f(\gamma z) = \chi(\gamma)j_\gamma (z)^k f(z) \]for all $\gamma \in \Gamma_0 (q)$ is called an \textit{automorphic function} of weight $k$, level $q$, and character $\chi$. The \textit{Laplace operator} of weight $k$ is defined by
\[ \Delta_k = y^2 \left( \frac{\partial^2}{\partial x^2} + \frac{\partial^2}{\partial y^2} \right) -iky\frac{\partial}{\partial x}, \] and a smooth automorphic function $f$ as above that is also an eigenfunction of the Laplace operator, i.e.,
\[ (\Delta_k + \lambda)f = 0 \]for some complex number $\lambda$, is called a \textit{Maass form} of corresponding weight, level, character, and the Laplace eigenvalue $\lambda$. One writes $\lambda(s) = s(1-s)$ and $s = 1/2 + ir$, with $r, s \in \mathbb{C}$, $r$ being known as the \textit{spectral parameter}. It is related to the Laplace eigenvalue $\lambda$ by the equation
\begin{equation}
   \lambda = \frac{1}{4} + r^2.  \label{eq:2.2.1}
\end{equation}

\subsection{Hecke eigenvalue}
\begin{definition}[~\cite{fg}]
The definition of the $n$th \textit{Hecke operator} $T_{n,\chi}, n \geq 1$ acting on the space of modular forms of level $q$, weight $k$, and character $\chi$ modulo $q$ is given by
\begin{equation}
   T_{n,\chi} : F(z) \longmapsto (T_{n,\chi}F)(z) = \frac{1}{n}\sum_{ad = n}\chi(a)a^k \sum_{b (\hspace{-0.2cm} \mod d)} F\left( \frac{az+b}{d} \right).   \label{eq:2.3.1}
\end{equation}For an eigenfunction $F$ of $T_{n.\chi}$, we shall denote the eigenvalue by $\lambda_F (n)$.
\end{definition}We denote the space of Maass forms of level $q$, weight $k$, and character $\chi$ modulo $q$ by $\mathcal{C}_k (q,\chi)$.
\begin{definition}[~\cite{fg}]
We define the action of the $n$th \textit{Hecke operator} $T^{\prime}_{n,\chi}$ on $\mathcal{C}_k (q,\chi)$ by
\begin{equation}
   T^{\prime}_{n,\chi} : f(z) \longmapsto (T^{\prime}_{n,\chi}f)(z) = \frac{1}{\sqrt{n}}\sum_{ad = n}\chi(a) \sum_{b ( \hspace{-0.2cm} \mod d)} f\left( \frac{az+b}{d} \right).    \label{eq:2.3.2}
\end{equation}
\end{definition}The operator \eqref{eq:2.3.2} is independent of $k$.\hspace{0.1cm}In Hecke theory for Maass forms, there is an orthonormal bases (Hecke basis) of Maass cusp forms consisting of forms that are common eigenfunctions of $T^{\prime}_{n,\chi}$ with $(n,q) = 1$. The forms in a Hecke basis is called \textit{Hecke-Maass cusp forms}. A Hecke-Maass cusp forms in the new subspace is called a \textit{newform} or a \textit{primitive from}. 


\subsection{The $L$-function attached to the twisted form ${\tf}$}
If $f$ is a Hecke-Maass form on $\text{SL}(2, \mathbb{Z})$ and $\chi$ is a primitive Dirichlet character modulo $q$, then the twisted form ${\tf}$ is a primitive Maass cusp form of level $q^2$. For normalized Hecke eigenvalues ${\la{f}{n}}$ of holomorphic forms $f,$ the $L$-function \textit{attached to the twisted form} ${\tf}$ is 
\begin{equation}
   L( s, {\tf} ) = \sum_{n} \frac{ \la{f}{n} \chi(n) }{n^s} = \prod_{p} L_{p} ( s,{\tf} )^{-1},   \label{eq:2.4.1}
\end{equation}where the local factor is
\begin{equation}
   L_{p} (s,{\tf}) = 1 - \la{f}{p} \chi(p)p^{-s} + \chi^2 (p)p^{-2s}.  \label{eq:2.4.2}
\end{equation}The Dirichlet series and the Euler product in \eqref{eq:2.4.1} are absolutely convergent for $\sigma > 1$.\hspace{0.1cm}We know from the theory of automorphic $L$-functions that the function $L(s,{\tf})$ has an analytic continuation to the whole complex plane and satisfies a functional equation relating the values at $s$ and $1-s$. 
\begin{theorem}[Theorem C in ~\cite{fg}]
There is an effectively computable absolute constant $c_2 > 0$ such that for any primitive form $f$ of some level $q$, spectral parameter $r$, and weight $k$, the $L$-function $L(s,f)$ does not vanish in the region
\begin{equation}
   \sigma \geq 1 - \frac{c_2}{\Ll}.      \label{eq:2.4.5}
\end{equation}
\end{theorem}Let $\Omega$ be the region in the complex plane given by
\begin{equation}
   \Omega = \left\{ \sigma + it \in \mathbb{C} : \sigma \geq 1 - \frac{c_2}{6\mathcal{L}} \right\},   \label{eq:2.4.3}
\end{equation}where
\begin{equation}
   {\Ll} = \log(q(|t| + |r| + 2)).   \label{eq:2.4.4}
\end{equation}We use these notations throughout  this paper.
\begin{lemma}[Lemma 4.1 in ~\cite{fg}]
Let $f$ be any Hecke-Maass cusp form for the full modular group and let $\chi$ modulo $q$ be any Dirichlet character. Let $c_2$ be the constant appearing in Theorem 2.4.1. Then, for every $s \in \Omega$, we have
\begin{equation}
   \frac{L^\prime (s,\tf)}{L(s,\tf)} \ll_f {\Ll},     \label{eq:2.4.6}
\end{equation}where the implied constant depends only on the form $f$.
\end{lemma}
\begin{lemma}[Lemma 4.3 in ~\cite{fg}]
Under the condition of Lemma 2.4.2, we have the uniform bound
\begin{equation}
   L(s,\tf) \hspace{0.3cm} \text{and} \hspace{0.3cm} L^{-1}(s,\tf)  \ll_f {\Ll}     \label{eq:2.4.7}
\end{equation}for all $s \in \Omega$, where the implied constant depends only on $f$.
\end{lemma}

\section{Proof of the main theorem}
\subsection{Preliminaries}
Let $f$ be the Hecke-Maass cusp form on $\text{SL}(2,\mathbb{Z})$ and $\chi$ be the primitive Dirichlet character modulo $q$. We consider the function $L(s,{\tf})$ as $F$ in \eqref{eq:2.1.1}. Then, \eqref{eq:2.1.2}-\eqref{eq:2.1.5} become
\begin{gather}
   \varphi(n,F) = \varphi(n,L) = n\prod_{p | n} \left( 1 - \frac{\alpha_1 (p)}{p} \right) \left( 1 - \frac{\alpha_2 (p)}{p} \right),   \label{eq:3.1.1}
   \\
   \gamma(p) = \la{f}{p}\chi(p) - \chi^2 (p)p^{-1},   \label{eq:3.1.2}
   \\
   C(F) = C(L) = \frac{1}{2}\prod_{p} \left( 1 - \frac{1}{p^2}( \la{f}{p}\chi(p) - \chi^2 (p)p^{-1} ) \right),   \label{eq:3.1.3}
\end{gather}and
\begin{equation}
   \alpha(n) = \mu(n)\prod_{p | n} (\la{f}{p}\chi(p) - \chi^2 (p)p^{-1})   \label{eq:3.1.4}
\end{equation}respectively. Here, the functions $\alpha_1$ and $\alpha_2$ in \eqref{eq:3.1.1} satisfy the following relations
\begin{gather}
   \alpha_1 (p) + \alpha_2 (p) = \la{f}{p} \chi(p),   \label{eq:3.1.5}
   \\
   \alpha_1 (p)\alpha_2 (p)  = \chi^2 (p).   \label{eq:3.1.6}
\end{gather}By Theorem 2.1.1, the asymptotic formula for $F(s) = L(s,{\tf})$ is
\begin{equation}
   \sum_{n \leq x} \varphi(n,L) = C(L)x^2 + E(x,L),   \label{eq:3.1.7}
\end{equation}where $E(x,L)$ is the associated error term. Using Theorem 2.1.1 again and partial summation, we can show the relation for $x \geq 1$,
\begin{equation}
   \sum_{n \leq x} \frac{\varphi(n,L)}{n} = 2C(L)x + E_0 (x,L),   \label{eq:3.1.8}
\end{equation}where $E_0 (x,L)$ is the associated error term.\hspace{0.1cm}Also, the function $H(s)$ in \eqref{eq:2.1.14} is determined for $\sigma > 1$ by the Dirichlet convolution as follows :
\begin{equation}
   H(s) = \sum_{n=1}^\infty \frac{1}{n^s} \sum_{d|n} \left( \mu(d)\prod_{p|d} \left( \la{f}{p} - \chi(p)p^{-1} \right)\chi(p) \right)\la{f}{\frac{n}{d}}\chi\left( \frac{n}{d} \right).  \label{eq:3.1.9}
\end{equation}

\subsection{Proof of Theorem 1.1.2}
Using Perron's formula and \eqref{eq:2.1.11}, we have  for $a > 1$
\begin{align}
    \sum_{n \leq x} \frac{\varphi(n,L)}{n} 
    &= \frac{1}{2\pi i}\int_{a-i\infty}^{a+i\infty} \left( \sum_{n=1}^\infty \frac{\varphi(n,L)}{n^{s+1}}  \right)\frac{x^s}{s}ds  \notag
    \\
    &=  \frac{1}{2\pi i}\int_{a-i\infty}^{a+i\infty} \zeta(s)\frac{H(s+1,f \otimes \chi)}{L(s+1, f \otimes \chi)} \frac{x^s}{s}ds.  \label{eq:3.2.1}
\end{align}Using Perron's formula again, we have for $a^\prime >1$
\begin{equation}
   x^{-1}\sum_{n \leq x} \varphi(n,L) =  \frac{1}{2\pi i}\int_{a^\prime - i\infty}^{a^\prime + i\infty} \zeta(s)\frac{H(s+1, f \otimes \chi)}{L(s+1, f \otimes \chi)} \frac{x^s}{s+1}ds.  \label{eq:3.2.2}
\end{equation}Combining \eqref{eq:3.2.1} with \eqref{eq:3.2.2}, we have
\begin{equation}
   \sum_{n \leq x} \left( 1 - \frac{n}{x} \right)\frac{\varphi(n,L)}{n} = \frac{1}{2\pi i}\int_{a-i\infty}^{a+i\infty} \zeta(s)\frac{H(s+1, f \otimes \chi)}{L(s+1, f \otimes \chi)} \frac{x^s}{s(s+1)}ds  \label{eq:3.2.3}
\end{equation}for $a > 1$. We calculate the integral on the right hand side of \eqref{eq:3.2.3}. Let $\mathscr{L}$ denote the contour 
\[ \mathscr{L} = \left\{ s = \sigma + it  : \sigma = -\frac{c_2}{6\mathcal{L}}, t \in \mathbb{R} \right\}, \]where the constant $c_2$  is the same as in Theorem 2.4.1. Since $\zeta(s)$ has a pole at $s=1$, the residue of the integrand in \eqref{eq:3.2.3} on the right hand side is
\begin{equation}
   \lim_{s \to 1} (s-1)\zeta(s)\frac{H(s+1, f \otimes \chi)}{L(s+1,f \otimes \chi)}\frac{x^s}{s(s+1)} = \frac{1}{2}\frac{H(2,f \otimes \chi)}{L(2,f \otimes \chi)}x.   \label{eq:3.2.4}
\end{equation}Since $H(s+1,f \otimes \chi)$ converges absolutely for $\sigma > -1/2$ by Theorem 2.1.3, $H(s+1,f\otimes \chi)$ has no poles in the considering region so that $H(s+1,f \otimes \chi)$ is bounded. Let $T$ be sufficiently large. By Lemma 2.4.3, the integral along ${\lL}$ is evaluated as 
\begin{align*}
   \left| \int_{{\lL}} \zeta(s)\frac{H(s+1, f \otimes \chi)}{L(s+1, f \otimes \chi)} \frac{x^s}{s(s+1)}ds \right|
   &\ll   q^2 \int_{-T}^T (\log(q(|t| + |r| +2)))^2 
   \\
   &\times (|t| + |r| + 2)^\frac{1}{2} \cdot |t|^{-2}\exp\left({-\frac{c_2}{6{\Ll}\log x}}\right)dt.
\end{align*}By taking $\log (q(|T| + |r| + 2)) = (\log x)^\frac{1}{2}$, it follows that $|T| = q^{-1}\exp\left({(\log x)^{1/2}}\right) \ll_{q,r} \exp\left({(\log x)^{1/2}}\right)$ and the above integral is bounded as
\begin{align*}
&\ll_{q,r} \int_{-T}^T (\log T)^2 \cdot T^\frac{1}{2} \cdot T^{-2} \cdot \exp\left( -\frac{c_2}{6}(\log x)^\frac{1}{2} \right)dt
\\
&= (\log T)^2 T^{-\frac{1}{2}} \exp\left( -\frac{c_2}{6}(\log x)^\frac{1}{2} \right)
\\
&\ll \exp\left( -\frac{c_2}{6}(\log x)^\frac{1}{2} \right)
\end{align*}for $T \geq 1$. Hence, 
\begin{equation}
    \int_{{\lL}} \zeta(s)\frac{H(s+1, f \otimes \chi)}{L(s+1, f \otimes \chi)} \frac{x^s}{s(s+1)}ds \ll  \exp\left( -\frac{c_2}{6}(\log x)^\frac{1}{2} \right).   \label{eq:3.2.5}
\end{equation}Also, the integral along the horizontal line from $-c_2 / 6{\Ll} +iT$ to $a+iT$ is 
\begin{equation} 
\leq   \int_{-c_2 / 6{\Ll} +iT}^{a+iT} \left| \zeta(s)\frac{H(s+1, f \otimes \chi)}{L(s+1, f \otimes \chi)} \frac{x^s}{s(s+1)} \right| |ds| \ll T^{-2}\log T.   \notag
\end{equation}The last bound tends to zero as $T \to \infty$. Hence the corresponding integral tends to zero as $T \to \infty$. Similarly, the integral along the opposite line tends to zero as $T \to \infty$. Therefore, by \eqref{eq:3.2.4},\eqref{eq:3.2.5} and the residue theorem, we have
\begin{align}
   \frac{1}{2\pi i}\int_{a-i\infty}^{a+i\infty} \zeta(s)\frac{H(s+1, f \otimes \chi)}{L(s+1, f \otimes \chi)} \frac{x^s}{s(s+1)}ds 
   &=  \frac{1}{2}\frac{H(2,f \otimes \chi)}{L(2,f \otimes \chi)}x   \notag
   \\
   &+ O\left( \exp\left( -\frac{ c_2 }{6}(\log x)^\frac{1}{2} \right) \right).  \label{eq:3.2.6}
\end{align}By \eqref{eq:3.1.7} and \eqref{eq:3.1.8}, we find from \eqref{eq:3.2.3} that
\begin{equation}
    \sum_{n \leq x} \left( 1 - \frac{n}{x} \right)\frac{\varphi(n,L)}{n} = C(L)x + E_0 (x,L) - x^{-1}E(x,L).  \label{eq:3.2.7}
\end{equation}Combining \eqref{eq:3.2.7} with \eqref{eq:3.2.6}, we have
\begin{align*}
   E_0 (x,L) 
   &= \left\{ \frac{1}{2}\frac{H(2,{\tf})}{L(2,{\tf})} - C(L) \right\}x + x^{-1}E(x,L)    \notag
   \\
   &+ O\left( \exp\left( -\frac{ c_2 }{6}(\log x)^\frac{1}{2} \right) \right). 
\end{align*}By \eqref{eq:2.1.13} and \eqref{eq:2.1.8}, we have
\begin{equation*}
\frac{H(2,{\tf})}{L(2,{\tf})} = \sum_{n=1}^\infty \frac{\alpha(n)}{n^2} = 2C(L)
\end{equation*}and we obtain the conclusion of Theorem 1.1.2 by taking $c_1 = c_2 / 6$.   $\square$

\begin{remark}
The reason why we take the function $L( s, {\tf} )$ whose Euler degree is equal to $2$ as $F$ in \eqref{eq:2.1.1} is because we know the estimate $L^{-1}( s, {\tf} )$ such as Lemma 2.4.3. If we know an estimate similar to Lemma 2.4.3, we will be able to obtain an asymptotic formula similar to Theorem 1.1.2 for the other $L$-functions whose Euler degree is equal to $2$.
\end{remark}

\textbf{Acknowledgments.}The author expresses his sincere gratitude to Prof. Kohji Matsumoto for giving me various advice in writing this paper.  

\end{document}